\documentclass[11pt,a4paper]{article}
\usepackage{amsmath,amssymb}

\usepackage{amsmath}
\usepackage{array,longtable}
\usepackage{color}
\usepackage[colorlinks,linkcolor=blue,citecolor=red]{hyperref}
\usepackage{graphicx}
\usepackage{wrapfig}

\DeclareMathSymbol{\bbbr}{\mathalpha}{AMSb}{"52}
\DeclareMathSymbol{\bbbc}{\mathalpha}{AMSb}{"52}

\newtheorem{theorem}{Theorem}

\newtheorem{definition}{Definition}

\newtheorem{lemma}[theorem]{Lemma}

\begin{document}

\title{Hexagonal 3-webs of circles with polar twisted cubic}

\author{{\Large Sergey I. Agafonov}\\
\\
Department of Mathematics,\\
S\~ao Paulo State University-UNESP,\\ S\~ao Jos\'e do Rio Preto, Brazil\\
e-mail: {\tt sergey.agafonov@gmail.com} }
\date{}
\maketitle
\unitlength=1mm

\vspace{1cm}

\begin{abstract}   
The paper describes and classifies  hexagonal circular 3-webs on unit sphere such that the polar points of the web circles lie on a twisted cubic, thus completing classification of hexagonal circular 3-webs with algebraic polar curves of degree three.

\bigskip

\medskip


\medskip

\end{abstract}



\section{Introduction} 

Consider three foliations in a plane domain such that the leaves of all foliations are circle arcs or straight line segments. If the foliations are transverse one to the other we say that these foliations form a circular 3-web in the domain. A 3-web is called hexagonal if at each point there is local diffeomorphism mapping the web in three families of parallel line segments. It is natural to study circular webs up to M\"obius transformations. Such transformations may map circles into straight lines therefore we treat straight line segments as arcs of circles of infinite radii. The difference between straight lines and circles disappears if we apply the inverse of stereographic projection to the web domain and obtain  a circular 3-web on the unit sphere.    

Stereographic projection puts the problem into a natural framework of the M\"obius geometry, instead of planar circular webs we study circular webs on the unit sphere $\mathbb{S}^2$.   Projective model of circle M\"obius geometry assigns points outside the unit sphere  to circles on this sphere: the  point is the polar point of the plane that cuts the circle on the sphere. 

Since we suppose that the web foliations are at least smooth each foliation gives a smooth arc outside the unit sphere.   Hence a circular 3-web on the unit sphere defines 3 curve arcs $\gamma_i$, $i=1,2,3$. 
   In what follows, the set of polar points of the web circles is called  the {\it polar curve} of the web. 
   
Given the polar curve, the web circles passing through a point  $p_0\in \mathbb{S}^2$ are cut on $\mathbb{S}^2$ by 3 planes $P_i$, polar to the points $p_i \in \gamma_i$ where the plane, tangent to $\mathbb{S}^2$ at $p_0$, meets the polar curve. The circles are transverse one to the other if the 3 lines $p_0p_i$ are transverse one to the other. Any parameter $u_i$ on the arc $\gamma_i$ is a first integral of the corresponding circular foliation. Thus the web is hexagonal if and only if the parameters can be chosen so that $u_1+u_2+u_3\equiv 0$, where $u_i$ is the value  giving $p_i \in \gamma_i$. The converse of the Abel Theorem (see \cite{BB-38}) immediately implies that, for hexagonal webs with planar polar curve, the arcs $\gamma_i$ are arcs of a planar cubic, possibly reducible.

The problem to describe hexagonal 3-webs formed by circles in the plane was posed by Blaschke and Bol in 1938 (see \cite{BB-38} p.31). By then Volk  \cite{V-29} and  Strubecker \cite{S-32}  have shown how to construct hexagonal circular 3-webs from  hexagonal linear  3-webs (i.e. webs formed by straight line segments). The construction, involving a central projection from a plane to a unit sphere followed by sterographic projection to a plane, was the first published  classification result, giving all hexagonal 3-webs with planar polar curves.   

 In fact, for all known examples the arcs $\gamma_i$ of the polar curve are algebraic. In the simplest case, when the arcs $\gamma_i$ are straight line segments, the web is formed by 3 pencils of circles.  The list of such webs was complete by the year 1977. The case of coplanar lines was settled by Volk and Strubecker,  the webs with non-coplanar lines appeared  in various publications since 1938, namely \cite{W-38,BB-38,B-73,E-74,L-77}. 
Finally, Shelekhov \cite{S-05} proved that the webs, described in these works, present a complete classification  of hexagonal webs whose polar curves splits into 3  non-coplanar lines.

In 1938, Wunderlich \cite{BB-38}  published a remarkable example of hexagonal circular 3-web. Its polar curve splits into 3 conics lying in 3 different planes. Since a generic plane, tangent to $\mathbb{S}^2$,  intersects 3 conics in 6 points, the Wunderlich web is actually 6-web, containing 8 hexagonal 3-subwebs. 

In 2014, Nilov \cite{N-14} found 5 new types of hexagonal circular 3-webs. For four of them, polar curves  split into a  line and a conic. The fifth is actually a 5-web whose polar curve is a union of a line and two conics. The line and two arcs on different conics give hexagonal 3-subwebs.  

The known examples suggest the problem to describe and to classify all hexagonal 3-webs with algebraic polar curves. It is natural to begin with polar curves of degree three. This first step is also motivated by the Graf and Sauer Theorem \cite{GS-24} which describes a linear hexagonal 3-web as being formed by tangents to a fixed curve of third class. In  the dual form, the construction of such webs is quite analogous to construction of circular webs from polar curves: take any cubic and consider the points $p_l$, dual  to the lines $l$, intersecting the cubic at 3 distinct  points $p_1,p_2,p_3$. Such points $p_l$ form a domain of a linear 3-web, the 3 lines of the web meeting at $p_l$ are dual to the points $p_1,p_2,p_3$. 

Volk  \cite{V-29} and  Strubecker \cite{S-32}) showed that planar polar curves of hexagonal 3-webs  are plane cubics. The case of three non-coplanar lines was solved by Shelekhov. All hexagonal 3-webs with polar curves, split into a smooth conic and a non-coplanar line, were described and classified by the author in \cite{A-25}: up to Möbius group action, there are 15 types, most of them depending on one parameter. (Four types of five in Nilov's paper \cite{N-14}  are, of course, webs  of this list.) 

The only case of a polar curve of degree three, which was not considered, is the one of a twisted cubic.  The example, presented in \cite{N-14} as hexagonal is actually not hexagonal (see Section \ref{prelim} for how to check this). 

The main result of the present paper solves this rest case and  completes classification of algebraic hexagonal circular 3-webs with polar curves of degree three, initiated in \cite{A-25}, where all 3-webs with reducible polar curves were listed. 

\begin{figure}[h]
\hspace{0.8cm}
\hspace{-0.5cm}\includegraphics[width=0.5\textwidth, trim={2.6cm 8cm 2cm 2cm}, clip]{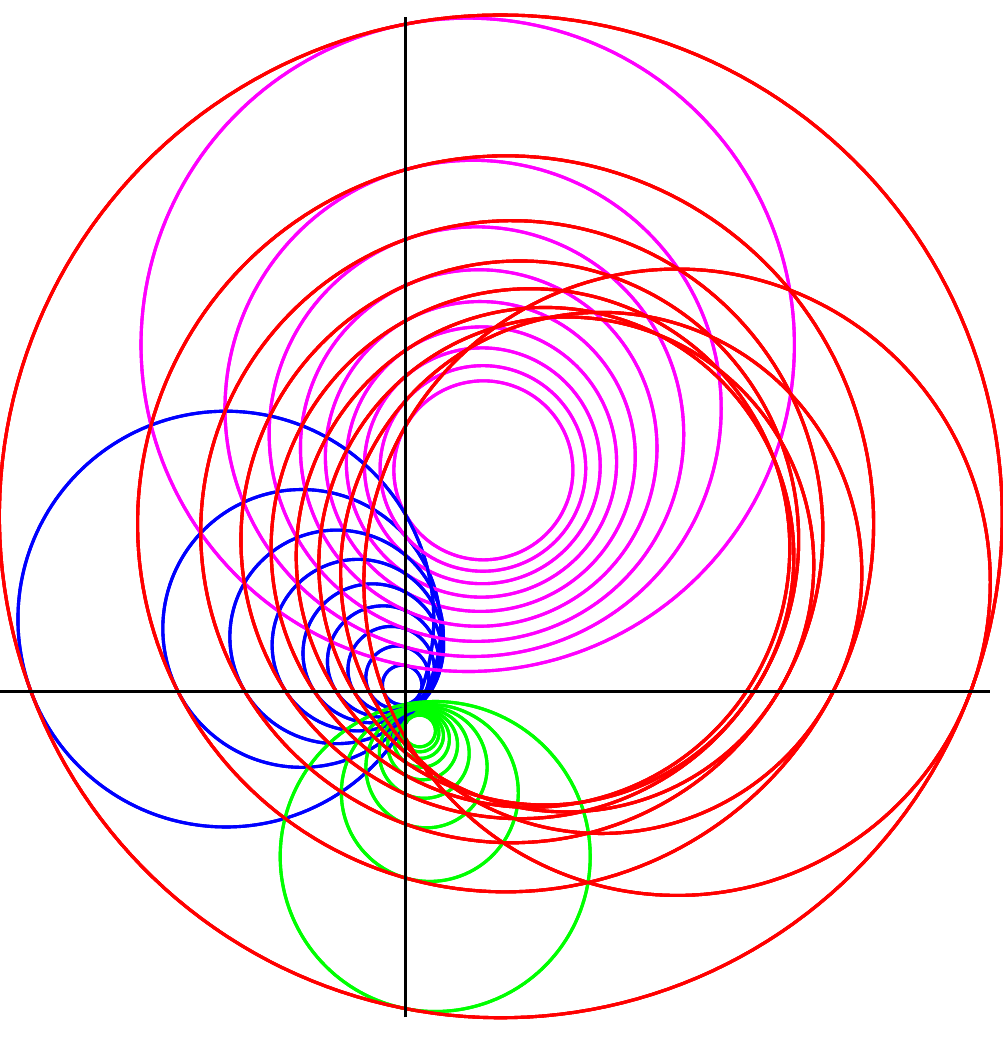} \includegraphics[width=0.54\textwidth, trim={2.7cm 8cm 0cm 2cm}, clip]{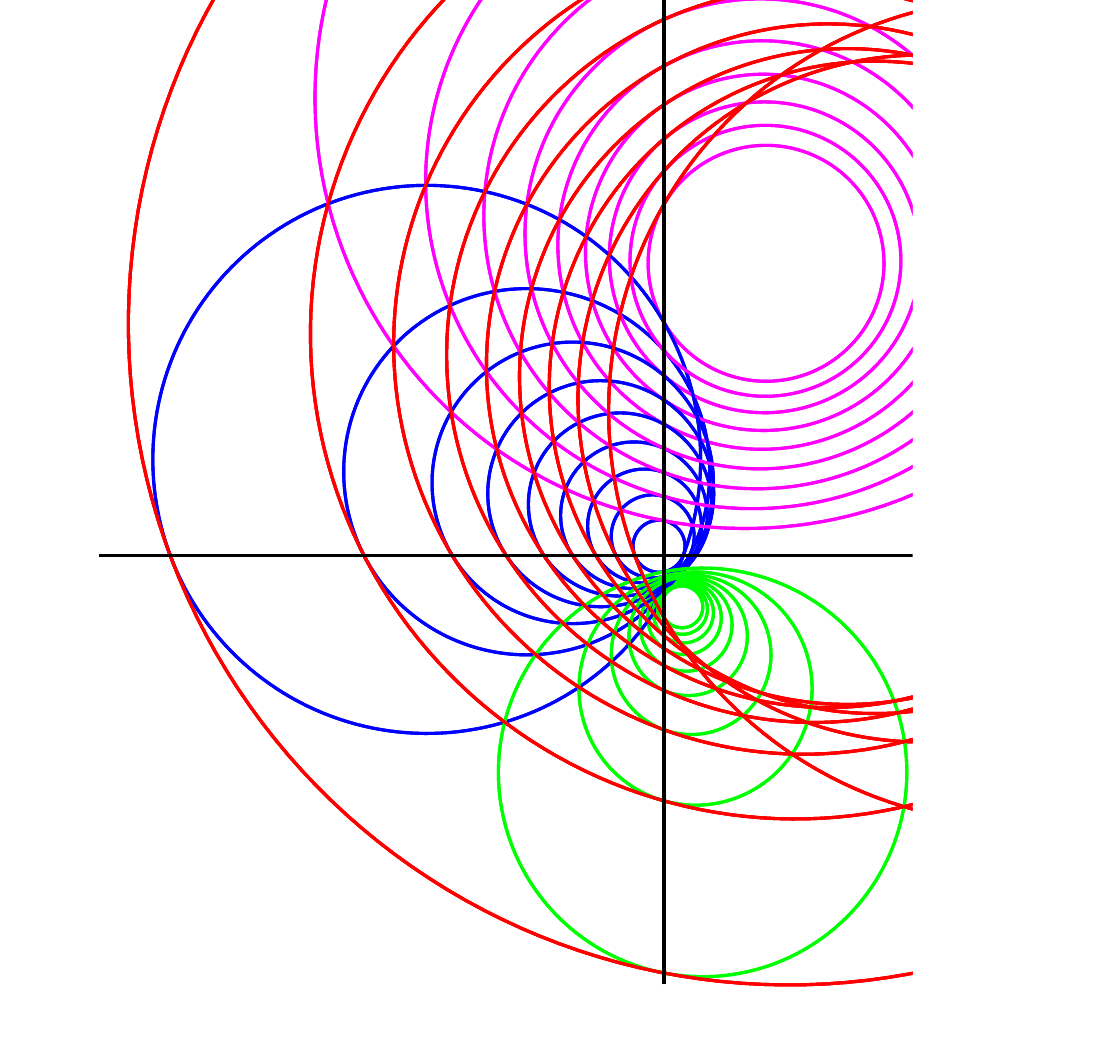}
 \caption{Stereographic image of 3-web with polar curve (\ref{cubic}), $x_0=\frac{\sqrt{3}}{2},  \ m=\frac{1}{\sqrt{3}}$.}\label{Fi0}
\end{figure}

\begin{figure}[h]
\hspace{0.8cm}
\includegraphics[width=0.85\textwidth, trim={0cm 8cm 0cm 2cm}, clip]{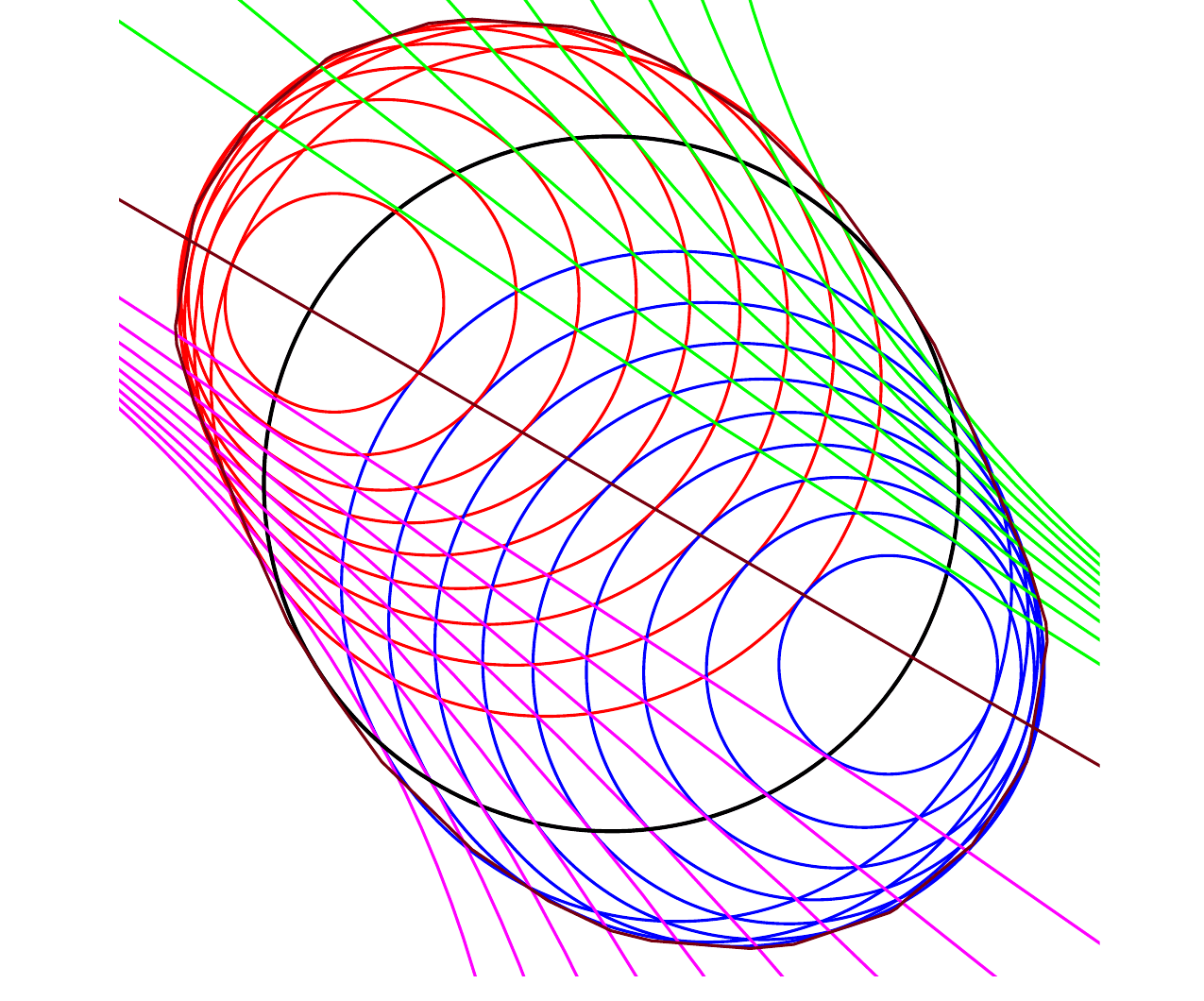}
 \caption{Stereographic image of 3-web with polar curve (\ref{cubic1}), $x_0=\frac{\sqrt{3}}{2}, \ y_0=\frac{1}{2}, \ m=\frac{1}{\sqrt{3}}$.}\label{Fi1}
\end{figure}

\begin{figure}[h]
\hspace{-0.6cm}
\includegraphics[width=1.0\textwidth, trim={0cm 8.5cm 0cm 3cm}, clip]{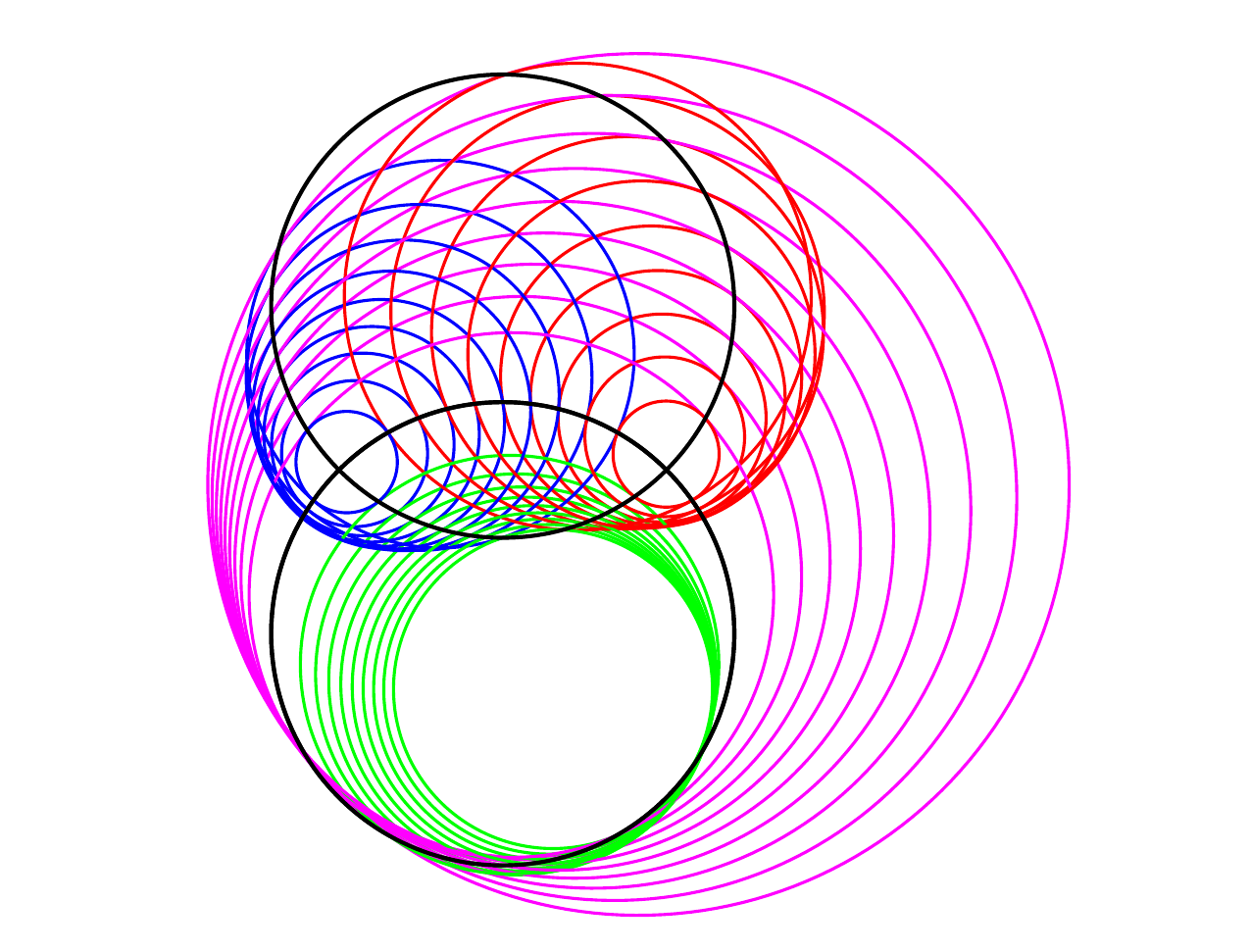}
 \caption{Stereographic image of 3-web with polar curve (\ref{cubic2}), $x_0=y_0=\frac{1}{\sqrt{2}}, \ m=\frac{1}{2}$.}\label{Fi2}
\end{figure}

\begin{theorem}\label{RNC}If the polar curve of a hexagonal circular 3-web is a twisted cubic then it is M\"obius equivalent to exactly one of the following  family, where $t$ is the parameter,  $m,x_0$ are positive real constants: 
\begin{equation}\label{cubic}
X=m(t^2-1),\ \ \ \ \
Y=t(t^2-1),\ \ \ \ \
Z=mx_0t,\ \ \ \ \
U=mx_0t^2.
\end{equation}
Generators of the homogeneous ideal for the twisted cubic (\ref{cubic}) can be chosen as 
$$
Z^2-U^2+x_0XU,\ \ \ \ \ 
  XU-mYZ, \ \ \ \ \ 
XZ+mY(x_0X-U).
$$
\end{theorem}
Stereographic projection image of a web from the family is presented in Figure \ref{Fi0}. 
Another two   representatives  are depicted in Figures \ref{Fi1} and  \ref{Fi2}, the M\"obius normalizations of these webs are different from (\ref{cubic}).   


The proof is based on exploring the simplest singular points of the web, namely, the points where exactly two circle arcs of different web foliations are tangent. 
Shelekhov observed \cite{S-05} that, for hexagonal 3-webs, such singular set is either a circle arc of the 3rd foliation or the common circle arc of the first two.   This restriction was rigid enough to isolate all the types in the classifications \cite{S-05} and \cite{A-25}. Moreover, it allowed to describe also all hexagonal circular 3-webs symmetric by one-parameter M\"obius subgroup (see \cite{A-25}). The polar curve of such symmetric webs are also algebraic, they split into lines and conics.

\section{ Web equation, singularities of hexagonal 3-webs, projective Möbius geometry}\label{prelim}

3-webs  in a two-dimensional domain are superpositions of 3 foliations $\mathcal{F}_i$ in general position. We will consider 3-webs on the unit sphere $\mathbb{S}^2$ cut by 3 families of planes.   
\begin{definition}
An algebraic circular web on the unit sphere $\mathbb{S}^2\in \mathbb{R}^3$ is a set of circles cut by planes dual  to points on some real (possibly reducible) algebraic curve $\gamma\in \mathbb{R}^3$, the duality being defined by the quadric $\mathbb{S}^2$. 

If one can choose 3 arcs $\gamma_i\in \gamma$ so that the three family of circles corresponding to points on these arcs form a 3-web $\mathcal{W}_3$ in some open domain $U\in \mathbb{S}^2$ we say that this 3-web $\mathcal{W}_3$ is algebraic circular 3-web and the union of irreducible components of $\gamma$ containing the arcs $\gamma_i$ is the polar curve of  $\mathcal{W}_3$.  

An algebraic circular 3-web is hexagonal if 
for any  point, were 3 circles from different foliations of the web are transverse one to the other,  there are a neighbourhood and a local diffeomorphism sending the web leaves of this neighbourhood  in 3 families of parallel line segments in the plane.   
\end{definition}
Topologically, hexagonality means an incidence property:  for any point $m$, each sufficiently small curvilinear triangle with the vertex $m$ and sides formed by the web leaves,  may be completed to the curvilinear hexagon, whose sides are web leaves and whose "large" diagonals are the web leaves meeting at $m$. One can choose 3 discrete families of leaves, one from each foliation, so that the chosen leaves form a triangulation of the domain, where they intersect transversally.  Hexagonality is equivalent to vanishing of the Blaschke curvature \cite{BB-38}. 

The hexagonality property of a circular  3-web can also be checked by the so-called {\it web equation} (see \cite{B-55}). It works as follows: let $u_1,u_2,u_3$ be local first integrals of the web foliations. Then these integrals are related by the {\it web equation}
\begin{equation}\label{webeq}
W(u_1,u_2,u_3)=0.
\end{equation}
At regular points, where all 3 leaves are transverse one to the other, this equation defines an implicit function $u_3=F(u_1,u_2)$. The web is hexagonal if and only if $F$ verifies the following partial differential equation of the 3rd order: 
\begin{equation}\label{hexPDE}
\partial_{u_1}\partial_{u_2}\ln\left(\frac{\partial_{u_1}F}{\partial_{u_2}F}\right)=0.
\end{equation}
In practice, calculations for algebraic circular 3-webs are performed via parametrization of the polar arcs $\gamma_1,\gamma_2,\gamma_3$ of the web. Namely, the parameters $u_1,u_2,u_3$ are first integrals, the equation (\ref{webeq}) manifests the fact that the plane containing the three points $p_i(u_i)\in \gamma_i$ is tangent to the unit sphere. This equation is obtained as follows. 
The stereographic projection maps the web circles into 3 families circles in the plane  
$$
\epsilon_i (x^2+y^2)+A_i(u_i) x+B_i(u_i) y+C_i(u_i)=0,
$$
(where $\epsilon_i=0$ or $\epsilon_i=1$, the case  $\epsilon_i=0$ giving a line segment). 
The corresponding polar points  with the tetracyclic coordinates
\begin{equation}\label{polarofcircle}
[A_i(u_i):B_i(u_i):C_i(u_i)-\epsilon_i:-C_i(u_i)-\epsilon_i]
\end{equation}
fill the arcs $\gamma_i$. Excluding $x,y$ between 3 equations (\ref{polarofcircle}) one gets the web equation. 

One can calculate all the partial derivatives in (\ref{hexPDE}) implicitly. 
 Then the web is hexagonal if and only if equation (\ref{hexPDE}) is true modulo (\ref{webeq}).       

Algebraic webs are global and they  inevitably have singularities as well as their 3-subwebs.  We call a singular  3-web hexagonal if the set of its {\it regular points} is non-empty and  its Blaschke curvature vanishes identically at regular points.  The point of a 3-web is regular if its 3 leaves are transverse at this point. 

\begin{lemma}\label{sing1}\cite{A-25}
Suppose that a local real or complex analytic hexagonal 3-web, formed by integral curves of three smooth direction fields $\xi_1,\xi_2,\xi_3$, has a singular point $p_0$ such that
\begin{enumerate}
\item all $\xi_i$ are well defined at $p_0$,
\item $\xi_1$ and $\xi_2$ are transverse at $p_0$,
\item $\xi_1=\xi_3$ at $p_0$, 
\end{enumerate}
then either the leaves of $\mathcal{F}_1$ and $\mathcal{F}_3$ trough $p_0$ coincide or 
$\xi_1=\xi_3$ along the leaf of  $\mathcal{F}_2$ through $p_0$. 
\end{lemma}
For algebraic circular 3-webs, the singularity at $p_0$,  described by the Lemma, means geometrically that  of 3 points $p_i$, where the plane $P_0$ tangent to $\mathbb{S}^2$ at $p_0$, meets the arcs $\gamma_i$, exactly two, namely $p_1$ and $p_3$, are collinear with $p_0$ and the plane $P_0$ does not contain any of the 3 lines tangent to $\gamma$ at $p_i$. The second condition means that $p_0$ is not on the envelope of any of the web foliations and therefore the algebraic direction fields $\xi_i$ remain smooth at $p_0$.      

Finally, we recall some necessary facts from hyperbolic 3D geometry. Following Blaschke \cite{B-29}, we call   the subgroup  of projective transformations of $\mathbb{RP}^3$, leaving invariant the quadric 
$$
X^2+Y^2+Z^2-U^2=0,
$$ 
the Möbius group. 

Infinitesimal generators of Möbius group  in homogeneous coordinates $[X:Y:Z:U]$ in $\mathbb{P}^3$, affine coordinates $x=\frac{X}{U},\ y=\frac{Y}{U},\ z=\frac{Z}{U}$ in $\mathbb{R}^3$ and cartesian coordinates $(\bar{x},\bar{y})$ in $\mathbb{R}^2$ related to points $(x,y,z)$ on the unit sphere via stereographic projection 
$$
x=\frac{2\bar{x}}{1+\bar{x}^2+\bar{y}^2},\ \ \ y=\frac{2\bar{y}}{1+\bar{x}^2+\bar{y}^2},\ \ \ z=\frac{1-\bar{x}^2-\bar{y}^2}{1+\bar{x}^2+\bar{y}^2},
$$
are as follows.
There are 3 rotations around the affine axes: 
$$
\begin{array}{l}
R_z=Y\partial_X-X\partial_Y=y\partial_x-x\partial_y=\bar{y}\partial_{\bar{x}}-\bar{x}\partial_{\bar{y}},\\
\\
R_x=Z\partial_Y-Y\partial_Z=z\partial_y-y\partial_z=\bar{x}\bar{y}\partial_{\bar{x}}+\frac{1}{2} (1-\bar{x}^2+\bar{y}^2)\partial_{\bar{y}},\\
\\
R_y=X\partial_Z-Z\partial_X=x\partial_z-z\partial_x=-\frac{1}{2} (1+\bar{x}^2-\bar{y}^2)\partial_{\bar{x}}-\bar{x}\bar{y}\partial_{\bar{y}},\\
\end{array}
$$
and 3 boosts (or "hyperbolic rotations"):
$$
\begin{array}{l}
B_x=U\partial_X+X\partial_U=\partial_x-x(x\partial_x+y\partial_y+z\partial_z)=\frac{1}{2} (1-\bar{x}^2+\bar{y}^2)\partial_{\bar{x}}-\bar{x}\bar{y}\partial_{\bar{y}},\\
\\
B_y=U\partial_Y+Y\partial_U=\partial_y-y(x\partial_x+y\partial_y+z\partial_z)= -\bar{x}\bar{y}\partial_{\bar{x}}+\frac{1}{2} (1+\bar{x}^2-\bar{y}^2)\partial_{\bar{y}},\\
\\
B_z=U\partial_Z+Z\partial_U=\partial_z-z(x\partial_x+y\partial_y+z\partial_z)=-\bar{x}\partial_{\bar{x}}-\bar{y}\partial_{\bar{y}}.\\
\end{array}
$$
The generators act on homogeneous coordinates as   the generators of $SO(3,1)$.

\section{Hexagonal circular 3-webs with twisted cubic polar  curves}

Let a circular hexagonal 3-web have a rational normal polar curve $N$. We complexify the real cubic $N$, the real  quadric $\mathbb{S}^2$ and consider the family of  bisecant lines $L$ of $N$ that are tangent to $\mathbb{S}^2$ at a point not on $N$. The  tangency points form a curve.  As the cubic $N$, being of degree three, cannot have trisecant lines, the closure $\gamma$ of the curve is a collection of circles $C_j$ (a circle in the complex context is the intersection of possibly complex plane, not tangent to $\mathbb{S}^2$,  with the complexification of $\mathbb{S}^2$).  The polar points $p_j$ of these circles lie on the complexification of $N$ by Lemma \ref{sing1}. 
\begin{definition}
The circles $C_j$ are called boundary circles. A point on a boundary circle is called singular for this circle if it is either on $N$ or the plane tangent to $\mathbb{S}^2$ at this point contains  the tangent line to $N$ at the polar point of $C_i$. Otherwise it is called regular point for $C_j$. 
\end{definition}
The choice of the  name {\it boundary circle} is due to the fact that a real domain of regular points of the web is bounded by arcs of these circles. Note that a point of boundary circle is singular for the circle if it does not satisfy the hypothesis of Lemma \ref{sing1}. 

Consider a boundary circle $C$, its polar point $p_0$ and the projection  $\pi: \mathbb{P}^3\setminus p_0 \rightarrow \ P_C$ from $p_0$ to the plane $P_C$ of the circle $C$. Then the closure of the image  of $N\setminus p_0$ is a conic, which we denote by $\pi(N)$, the intersection $p'_0$ of the tangent line to $N$ at $p_0$ with $P_C$ is on $\pi(N)$ and the osculating plane of $N$ at $p_0$ is tangent to $\pi(N)$ at $p'_0$.     

\begin{lemma}\label{projection} If $N$ is a twisted cubic polar curve of hexagonal 3-web  then 
\begin{enumerate}
\item if $\pi(N)$ and $C$ are tangent at $q$ then $q=p'_0$,   
\item if $\pi(N)$ intersects $C$ transversally at a point $q$ then the line $p_0q	$ is a bisecant of $N$ and the tangent line to $C$ at $q$ contains $p'_0$
\item if $p'_0\notin C$ then the points, where the tangents from $p'_0$ to $C$ touch $C$, are on $\pi(N)$.    

\end{enumerate} 
\end{lemma}
{\it Proof:} 
A generic  line $L$ of the considered family of bisecants of $N$ touches $\mathbb{S}^2$ at a regular point $p_L$ for the circle $C$  and meets the curve $N$ at two points $p_1(L),p_2(L)$. The image $\pi(L)$ of $L$ is a tangent line  to $C$ at $p_L$. 
Consider the family of planes $P_L$  containing $p_0$ and $L$.  

Let $\pi(N)$ and $C$ be tangent at $q$.  As $p_L$ approaches its limit position $q$, the points $\pi(p_1(L)),$ $\pi(p_2(L))$ also tend to $q$. As $N$ does not have trisecant lines (tangent lines to $N$ intersecting $N$ an some other point count as trisecants), both points $p_1(L),p_2(L)$ tend to  
$p_0$, the plane $P_L$ to the osculating plane of $N$ at $p_0$ and $p_L$ to $p'_0$. Hence $q=p'_0$. 

 If $\pi(N)$ intersects $C$ transversally at a point  $q\notin N$ we consider again the family of planes $P_L$. As $p_L$ approaches its limit position $q$, one of the points $\pi(p_1(L)),$ $\pi(p_2(L))$, say $\pi(p_1(L))$, tends to $q$ and the other to some point $r$ on $\pi(N)$. If $r\ne p'_0, \ q\ne p'_0$ then the lines $p_0r$ and $p_0q$ are bisecant lines of $N$ and meet it at some points $p_1\ne p_0$, $p_2\ne p_0$. On the other hand, the points $p_1,p_2$ must be collinear with $q$ which is not possible. 
 
 If $q=p'_0$ then the plane $P_q$ is tangent to $N$ and  the intersection of $C$ and $\pi(N)$ is not transversal at $q$. Hence $r=p'_0$ and the plane $P_q$, being not osculating of $N$, meets $N$ at some point on the bisecant line $p_0q$ of $N$.   In fact, as $P_L$ approaches the position of a plane tangent to $N$ at $p_0$ and $p_1(L)$ tends to $p_0$, the limit position of the point $p_2(L)$ is collinear with $q$ and $p_0$.

Now suppose than $p'_0\notin C$. 
As $P_L$ approaches the position of a plane tangent to $N$ at $p_0$, one of two points $p_1(L),p_2(L)$, say $p_1(L)$ tends to $p_0$. 
The limit position of $p_2(L)$ is not $p_0$, therefore, being collinear with $p_L$ and $p_1(L)$, the point $p_2(L)$ must tend to some point on a line joining $p_0$ with the limit  position of $p_L$. Thus, in the limit, the line $p_0p_2(L)$ is bisecant of $N$ intersecting $C$. The Lemma is proved.      
\hfill $\Box$
\begin{lemma}\label{Poncelet}
Let $N$ be  polar twisted cubic curve of a hexagonal 3-web, $C$ its boundary circle with polar point $p_0$, and $\pi(N)$ the projection of $N$ from $p_0$ to the plane $P_C$ of $C$. If there are two distinct points $q_1,q_2$ in the intersection $N\cap P_C$ that are not on $C$ then the third point $q_3$ in this intersection is also not contained in $C$ and the circle $C$ is inscribed in the triangle $q_1q_2q_3$.  
\end{lemma}
{\it Proof:} At least one of the points $q_1,q_2$ is different from $p'_0$, let it be $q_1$. Then the tangent lines from $q_1$ to $C$ meets $\pi(N)$ at two points on $N$. Hence one of them is $q_2$ and the other some point $q_3$. Repeating this reasoning for $q_2,q_3$ we deduce the desired conclusion.   
\hfill $\Box$

\begin{lemma}\label{circles} For polar twisted cubic curve of hexagonal 3-web  the following facts are true.
\begin{enumerate} 
\item The curve $\gamma$ consists of 2 circles $C,\bar{C}$. 
\item The  circles $C,\bar{C}$ intersect in two distinct points $c_1,c_2$, lying on $N$. 
\item The polar points $p_0,\bar{p}_0$ of $C,\bar{C}$ are conjugated with respect to $\mathbb{S}^2$. 
\end{enumerate} 
\end{lemma}
{\it Proof:} 
The curve $\gamma$ cannot contain more than 2 circles. To show this, consider three boundary circles $C_1,C_2,C_3$. As each of them contains at most 5 singular points there is a circle $c$ of the web, intersecting them in 6 distinct regular  points $r_i$. Note that there is no point common for all web circles with polar points on $N$ as $N$ is not planar, therefore $c$ can be chosen so that it misses the points where a pairs of $C_1,C_2,C_3$ intersect.   Consider the cone of lines joining the polar  point $p_c$ of $c$ with the points of $c$. Each of 6 lines $p_cr_i$ of the cone meet $N$, which is not possible: these lines would be in the intersection of the cone with another quadratic cone with the vertex $p_c$, namely, the one over the projection of $N$ from $p_c$.   In what follows, we call this pair of cones  a {\it pair of cones at $p_c\in N$} for brevity. 
  
 If $\gamma$ consists of one circle $C$ then we choose again a circle $c$ of the web, intersecting $C$ at 2 regular points $r_1,r_2$ and consider the pair of cones at $p_c$. The cones intersect at two lines $p_cr_1$ and $p_cr_2$ therefore at least one, say $p_cr_1$ is multiple. Since it is true for an open set of points $p_c\in N$ it holds also for the polar point $p_0$ of $C$. Thus the conics, where the cones intersect the plane $P_C$ of the circle $C$, intersect $C$ at most at two points $r_1,r_2$ with exactly one point, say $r_1$, being multiple. Due to the claim 1 of Lemma \ref{projection}, the multiple point is $p'_0$. 
 
Suppose that the projection $\pi(N)$ from $p_0$ to the plane $P_C$ of the circle $C$ intersects $C$ in just one point $r$ of multiplicity $4$. Using a  suitable M\"obius transform, we move $p_0$ to $p_z=[0:0:1:0]$, the point $r$ to $[1:0:0:1]$. Then in the affine coordinates $x,y$ in the plane $P_C$, given by $Z=0$, the conic $\pi(N)$ has the form $y^2+2x-2+\lambda (x-1)^2=0$. Applying a suitable complexified  M\"obius transform, preserving $C$, we can move the intersection point of this conic with the $x$-axis to infinity, thus bringing the conic to the form $y^2+2x-2=0$. 
Consider the points of intersection $P_C\cap N$. All of them are on $\pi(N)$ and none of them coincides with $r=p'_0$. By the Poncelet porism,  there is no triangle, inscribed in $\pi(N)$ and circumscribed about the circle $C$: it is enough to try the one with one side touching $C$ at the point $[-1:0:0:1]$. Therefore all the points in $P_C\cap N$ coincide by Lemma \ref{Poncelet}. Let us denote this triple point by $q$. Considering again the family of planes $P_L$ approaching one of two points where the tangent lines from $q$ to $C$ touches $C$, we conclude that the corresponding tangent is also a tangent to $N$ at $q$. Thus we get a contradiction, since there are two tangents from $q$ to $C$.      
   
Suppose that the projection $\pi(N)$ from $p_0$ to the plane $P_C$ of the circle $C$ intersects $C$ in two points,  $r_1$ of multiplicity $3$ and $r_2$ of multiplicity 1. Then, by Lemma \ref{projection}, the tangent line to $C$ at $r_2$ contains $p'_0=r_1\in C$ which is not possible.

Thus $\gamma$ consists  of exactly two circles $C,\bar{C}$ with polar points $p_0,\bar{p}_0$. 
If the circles $C,\bar{C}$ are tangent at $c$, then the points $c,p_0,\bar{p}_0$ are collinear and  $c\notin N$. The plane tangent to $\mathbb{S}^2$ at $c$ contains a third point $q$ of its intersection with $N$. If the point $q$ does not coincide with none of $p_0,\bar{p}_0$ then the circle with polar point $q$ would be the third boundary circle which is not possible. Therefore one of the points $p_0,\bar{p}_0$ is a multiple intersection point. Let it be $\bar{p}_0$. Now the projection $\pi(N)$ from $p_0$ on the plane of circle $C$ is tangent to $C$ at $\pi(\bar{p}_0)$.  
By the claim 1  of Lemma \ref{projection}, the point $p'_0$ coincides with $c$, the tangent line to $N$ at $p_0$ contains $\bar{p}_0$ and we end with a trisecant line of $N$, the point $p_0$ counted twice. 

Thus the boundary circles $C$ and $\bar{C}$ intersect in two distinct  points $c_1,c_2$.  None of them can be regular for both circles: otherwise the third point of intersection of the plane $P_i$, tangent to $\mathbb{S}^2$ at the point $c_i$, regular for both boundary circles, with $N$ would lie on both lines $p_0c_i$ and $\bar{p}_0c_i$ crossing at $c_i$.  If $c_1\notin N$ then the plane $P_1$ must be tangent to $N$ at exactly  one of points $p_0,\bar{p}_0$. Let it be $\bar{p}_0$. Then $c_1$ is regular for $C$ and as the plane $P_L$ tends to $P_1$ the points $p_1(L)$, $p_2(L)$ tend to $\bar{p}_0$ and the line $p_1(L)p_2(L)$ to the tangent to $N$ at $\bar{p}_0$. Thus the tangent to $N$ at $\bar{p}_0$ contains $c_1$, therefore the projection $\bar{\pi}(N)$ of $N$ from $\bar{p}_0$ to the plane of $\bar{C}$ is tangent to $\bar{C}$ which is not possible as the tangent line to $\bar{C}$ at $c_1$ contains the 3rd point on $\bar{\pi}(N)$, namely, the  image $\bar{\pi}(p_0)$ of $p_0$.   

Hence $c_1,c_2$ lie on $N$.   The plane $P_C$ of $C$ meets $N$ at $c_1$ and $c_2$. If the third point $q$ in $P_C\cap N$ is distinct from $c_1$ and $c_2$ then the tangent lines from $q$ to $C$ must meet $N$. Thus these  lines are tangent to $C$ at $c_1$ and $c_2$ and $q$ is contained in the line $P_1\cap P_2$. But this line contains 2 points of $N$, namely $p_0$ and $\bar{p}_0$. Hence $q=\bar{p}_0$ and $p_0,\bar{p}_0$ are conjugated.  If the third point $q$ coincides with one of $c_1,c_2$, say $c_1$, then the plane $P_C$ is tangent to $N$ at $c_1$ and as the plane $P_L$ tends to $P_1$, the points $p_1(L)$, $p_2(L)$ tend to $\bar{p}_0$ and $c_1$.  At the limit position, the line $p_2(L)p(L)$ coincides with the line $\bar{p}_0c_1$. In this process, the plane $p_2(L)p(L)c_1$ tends to $P_C$.
 Thus $\bar{p}_0\in P_C$ and $p_0,\bar{p}_0$ are conjugated. 
\hfill $\Box$

\vspace{7pt}

The necessary restrictions imposed by  Lemmata \ref{sing1},\ref{projection},\ref{Poncelet},\ref{circles} turn out to be also sufficient for hexagonality.

\vspace{7pt}

{\noindent \it \bf Proof of Theorem \ref{RNC}:}  First consider the case when the bisecants $L$ of $N$, tangent to $\mathbb{S}^2$, are real and therefore the two boundary circles $C,\bar{C}$ are also real. 
Applying a suitable M\"obius transformation, we send  $p_0$ to $p_y=[0:1:0:0:0]$,  $\bar{p}_0$ to $p_x=[1:0:0:0]$, $c_1$ and $c_2$ to $[0:0:\pm 1:1]$. Then the projection $\pi(N)$ from $p_0$ contains the points $[0:0:\pm 1:1]$ and $p_x$. Therefore the point $p'_0$, lying also on  $\pi(N)$,  is not contained in neither of lines $z=\pm 1$ in the plane $Y=0$: otherwise the conic $\pi(N)$ would intersect one of these lines in 3 points. If the affine coordinate $z$ of  $p'_0$ is inside the interval $(-1,1)$ then, applying a transform from the M\"obius subgroup generated by $B_z$, one can send $p'_0$ to some point on the $x$-axes, otherwise one can send $p'_0$ to some point at the infinite line of the plane $Y=0$. The subgroup preserves the chosen positions of $p_0,\bar{p}_0,c_1,c_2$. 

Consider first the case when $p'_0$ is on the $x$-axes. Note that it cannot coincide with $p_x$ since $N$ cannot have trisecant lines. Then $p'_0=[1:0:0:x_0]$, where one may suppose $x_0>0$ without loss of generality.   
The conic $\pi(N)$ is uniquely defined by its 5 points $p_x$, $[0:0:\pm 1:1]$ and  $[x_0:0:\pm z_0:1]$, the last two being the points where tangents from $p'_0$ meet $C$. (Here $z_0$ is complex for $x_0>1$).  In the homogeneous coordinates $X,Z,U$ on the plane  $Y=0$, the conic has the form
\begin{equation}\label{con1}
Z^2+x_0XU=U^2.
\end{equation}
We choose the  parameter $t$ of a rational polynomial in $t$ parametrization  $\Gamma: \mathbb{RP}^1\to \mathbb{RP}^3$, $t\mapsto [X(t):Y(t):Z(t):U(t)]$ of $N$ so that $\Gamma(\infty)=p_y$, $\Gamma(\pm 1)=[0:0:\pm 1:1]$. Then $Y(t)$ is a polynomial of degree 3 and  $X(t),Z(t),U(t)$ are at most quadratic in $t$. Moreover, $X(t),Y(t)$ have two roots  $t_{\pm}=\pm 1$ and therefore  $X(t)=m(t^2-1)$, where $m\ne 0$ can be chosen positive, $Y(t)$ is also divisible by $t^2-1$, and $Z(t)=P(t^2-1)+Kt+L$, $U(t)=Q(t^2-1)+Lt+K$. Now the condition $p'_0=[1:0:0:x_0]$ gives $P=0$ and $Q=mx_0$. The curve $N$ contains $p_x$, i.e.  $\Gamma(t_0)=[1:0:0:0]$. Hence $Y=(t^2-1)(t-t_0)$, $K=Q=mx_0$, $L=-Kt_0=-mx_0t_0$. The projection $\pi(N)$ is the conic (\ref{con1}). This gives $t_0=0$ and the parametrization assumes the form (\ref{cubic}).  
One checks by direct computation that the corresponding web is hexagonal.  

The projection $\bar{\pi}(N)$ from $\bar{p}_0$ to the plane $X=0$ of the circle $\bar{C}$ has he form 
\begin{equation}\label{con2}
Z^2+mx_0YZ=U^2,
\end{equation}
with $\bar{p}'_0=[0:-1:mx_0:0]$ at infinity. We show that the second, not yet considered case with $p'_0$ at infinity, is  another "face" of the found web. One has just to swap $p_0$ and $\bar{p}_0$. 
To this end we start with the M\"obius normalization with  $c_1,c_2=[0:0:\pm 1:1]$, $p_0=p_x$, $\bar{p}_0=p_y$. 
The conic $\pi(N)$ is uniquely defined by its 5 points $p_y$, $[0:0:\pm 1:1]$ and  $[0:y_0:z_0:1]$, $[0:-y_0:-z_0:1]$,  the last two being the points where tangents from $p'_0=[0:z_0:-y_0:0]$ at infinity meets $C$. Note that $z_0\ne 0$ and, without loss of generality, one may suppose that $y_0,z_0$ are positive.   In the homogeneous coordinates $Y,Z,U$ on the plane  $X=0$, the conic has the form
\begin{equation}\label{con3}
Z^2+\frac{y_0}{z_0}YZ=U^2.
\end{equation}
We choose the  parameter $t$ of a rational parametrization  $\Gamma: \mathbb{RP}^1\to \mathbb{RP}^3$, $t\mapsto [X(t):Y(t):Z(t):U(t)]$ of $N$ so that $\Gamma(0)=p_0=p_x$, $\Gamma(\pm 1)=[0:0:\pm 1:1]$. Then   $Y(t),Z(t),U(t)$ have the factor $t$. Moreover, $X(t),Y(t)$ have two roots  $t_{\pm}=\pm 1$ and therefore  $X(t)=(t^2-1)(a+bt)$, where $a\ne 0$ can be chosen positive, $Y(t)=t(t^2-1)$, and $Z(t)=Bt+(P+Q)t^2+Dt^3$, $U(t)=Pt+(B+D)t^2+Qt^3$. Since $p_y=[0:1:0:0]\in N$ and the corresponding value $t_y$ of $t$ is not $\pm 1$ we conclude that $t_y=\infty$ and therefore $D=Q=0$.  
Now the condition $p'_0=[0:z_0:-y_0:0]$ gives $P=0$, $B=y_0/z_0$ and the parametrization assumes the form (\ref{cubic}) with $a=m$ and $y_0/z_0=mx_0$. Observe that $y_0$ and $z_0$ are subject to $y^2_0+z^2_0=1$ and therefore the value $y_0/z_0$ defines both of them.

To separate the orbits of the M\"obius action on the polar cubics, we consider  two discrete invariants 
$$
S={\rm sign}(c_1,c_2)(c_2,p_0')(p_0',c_1), \ \ \ \ \ \bar{S}={\rm sign}(c_1,c_2)(c_2,\bar{p}_0')(\bar{p}_0',c_1),
$$ 
and two continuous invariants
$$
I=\frac{(p_0,p_0)(\bar{p}'_0,\bar{p}'_0)}{(p_0,\bar{p}'_0)^2}, \ \ \ \ \ \bar{I}=\frac{(\bar{p}_0,\bar{p}_0)(p'_0,p'_0)}{(\bar{p}_0,p'_0)^2},
$$
where $(\xi,\eta)$ is the scalar product with signature (3,1), defining $\mathbb{S}^2$ as a quadric in $\mathbb{P}^3$, and the involved points  $\xi,\eta\in \mathbb{P}^3$ represent vectors in $\mathbb{R}^4$ by their homogeneous coordinates (see \cite{B-29} for more detail).  
Fore the cubic (\ref{cubic}), one has 
$$
S=-1,\ \ \ \bar{S}=+1,\ \ \ I=1+m^2x_0^2, \ \ \ \bar{I}=1-x_0^2.
$$ 
Thus $m^2$ and $x_0^2$ are also M\"obius invariants effectively separating the orbits of cubics (\ref{cubic}) since $x_0,m$ are both positive. Moreover, the discrete invariants $S,\bar{S}$ distinguish the two "faces" of the cubic, "visible" from $p_0$ and $\bar{p}_0$.

Now we prove that the case  of complex bisecant lines of $N$ tangent to $\mathbb{S}^2$ is not possible for real hexagonal 3-webs. 
The key observation is that the two families of complex bisecants $L$ of the complexification of a real cubic  $N$, with $L$ also tangent to the complexifications of the real sphere $\mathbb{S}^2$ along two circles $C$ and $\bar{C}$, are complex conjugated. 

Consider a real hexagonal 3-web $\mathcal{W}_3$ with real polar twisted cubic $N$,  whose bisecants, tangent to $\mathbb{S}^2$, are complex. The points $p_0,\bar{p}_0$ are complex conjugated, as well as the corresponding polar planes. Being conjugated, these planes intersect in a real line, containing the complex intersection points of the complex circles  $C$ and $\bar{C}$. Therefore this line is elliptic and its dual, containing complex points $p_0,\bar{p}_0$, is hyperbolic.  We can move this hyperbolic line to the $z$-axis $X=Y=0$. Then $p_0=[0:0:z_0:1]$ and $\bar{p}_0=[0:0:z_0^*:1]$, where $z_0^*$ is the complex conjugated of $z_0$, moreover, holds $z_0z_0^*=1$. One easily checks that there is a M\"obius transform, namely a hyperbolic rotation,  moving $[0:0:z_0:1]$ to $[0:0:i:1]$ where  $i^2=-1$. 
The line dual to the $z$-axis cuts $\mathbb{S}^2$ at points $[\pm i:1:0:0]$. 
Hence $c_1=[ i:1:0:0]$ and $c_2=[- i:1: 0:0]$. Consider a real rational parametrization $\Gamma: \mathbb{RP}^1\to \mathbb{RP}^3$, $t\mapsto [X(t):Y(t):Z(t):U(t)]$ of $N$. Since $c_1,c_2$ are complex, the corresponding parameter values are also complex and conjugated. By a suitable real linear-fractional transform we can send them to $\pm i$. Conditions $\Gamma(i)=c_1,$ $\Gamma(-i)=c_2$ imply that the polynomials $Z(t)$ $U(t)$ are divisible by $t^2+1$.  Similarly, if $\Gamma(\alpha\pm i\beta)=[0:0:\pm i:1]$, we can use a change of parameter $t$, preserving $i$, to move $\alpha\pm i\beta$ to $\pm ip$, where $p$ is real and nonzero. Thus $X(t),Y(t)$ have the factor $t^2+p^2$. Taking into account the other two components of $\Gamma$, we get 
$$
\begin{array}{c}
X(t) = (t^2+p^2)(Kt - L), \ \ \ Y(t) = (t^2+p^2)(Lt + K), \\
\\
Z(t)=(t^2 + 1)(Mt - pN), \ \ \ U(t)=(t^2 + 1)(Nt+pM)
\end{array} 
$$ 
 with real $K,L,M,N$. 
 
 Now we use restrictions imposed by Lemma \ref{projection} and consider the projection of $N$ from $p_0$ on its dual plane $U=iZ$. Let us choose $X,Y$ and $U$ as the homogeneous coordinates on this plane.  In these coordinates, the projection $\pi(N)$ is parametrized by 3 quadratic polynomials as follows $t\mapsto \left[\frac{X(t)}{it+p}:\frac{Y(t)}{it+p}:\frac{iZ(t)+U(t)}{it+p}\right]$, and the circle $C$ is $X^2+Y^2-2U^2$. 
One computes 
$$
p'_0=[2ip(pKi - L):2ip(pLi +K):\frac{1}{2}(p^2-1)(Ni-M):\frac{i}{2}(p^2-1)(Ni-M) ]
$$ 
and checks that the restriction of item 3 of Lemma \ref{projection} implies that  $K^2+L^2$ vanishes. Since $K,L$ are real, both must vanish. But then the map $\Gamma$ parametrizes the line $X=Y=0$ and not a twisted cubic. 
\hfill $\Box$\\

 \section{Concluding remarks} 

The M\"obius normalization (\ref{cubic}) is not the best choice for the visualization, as presented in Figure \ref{Fi0}. Better pictures were obtained for the curve  

\begin{equation}\label{cubic1}
\begin{array}{l}
X=my_0\left(t^2 + (x_0^2 - y_0^2)/4\right),\ \ \ \
Y=-mx_0\left(t^2 +(y_0^2 - x_0^2)/4\right),\\
\\
Z=t^3-t/4,\ \ \ \ 
U=-mx_0y_0t,\ \ \ \  x_0>y_0,

\end{array}
\end{equation}
in  Figure \ref{Fi1},  
 or the curve 
\begin{equation}\label{cubic2}
 \begin{array}{l}
X=x_0mt,\ \ \ \, 
Y=m(-x_0^2t^2+y_0^2+1)/2y_0,\\
\\
Z=(t-t^3)/8,\ \ \ \
U=m(x_0^2t^2+y_0^2+1)/2,
\end{array}
\end{equation}
in  Figure \ref{Fi2}.   
In both parametrizations, $x_0,y_0,m$ are positive  subject to  $x_0^2+y_0^2=1$.
Generators of the homogeneous ideal for the twisted cubic (\ref{cubic1}) can be chosen as 
$$
\frac{X^2}{y_0^2}-\frac{Y^2}{x_0^2}+\left(\frac{1}{x_0^2}-\frac{1}{y_0^2}\right)U^2,\ \
  2mZ \left(\frac{X}{y_0}+\frac{Y}{x_0}\right)-U\left(\frac{X}{x_0}+\frac{Y}{y_0}\right), \ \
 \left(\frac{X}{y_0}-\frac{Y}{x_0}\right)^2+\frac{4m}{x_0y_0}ZU-\frac{1}{x_0^2y_0^2}U^2. 
$$
For the cubic (\ref{cubic2}), homogeneous generators of the ideal are
$$
(1+y_0^2)X^2+y_0^2Y^2-U^2,\ \
\frac{2y_0^2}{mx_0^3}XU+8ZU+8y_0YZ-\frac{2y_0}{mx_0^3}XY,\ \
x_0^2X^2-y_0^2Y^2+2y_0YU-8mx_0^3XZ-U^2.
$$

\section*{Acknowledgements}
This research was supported by FAPESP grant \# 2022/12813-5.


\begin{thebibliography}{99}


\bibitem[A-25]{A-25} Agafonov, S. I.  Hexagonal circular 3-webs with reducible polar curves of degree three, SIGMA Symmetry Integrability Geom. Methods Appl. 21 (2025), 043, 31p, 
 doi.org/10.3842/SIGMA.2025.043


\bibitem[B-73]{B-73}  Balabanova, R. S. Hexagonal circular three-webs,  Plovdiv. Univ. Nauchn. Trud. Mat.,  pp. 128-141, 11, No. 4, 1973.

\bibitem[B-29]{B-29} Blaschke, W. 
Vorlesungen über Differentialgeometrie und geometrische Grundlagen von Einsteins Relativitätstheorie. III: Differentialgeometrie der Kreise und Kugeln. Bearbeitet von G. Thomsen. (German)  Berlin, J. Springer (Grundlehren der mathematischen Wissenschaften in Einzeldarstellungen) (1929).

\bibitem[BB-38]{BB-38} Blaschke, W.,  Bol, G.
 Geometrie der Gewebe, Topologische Fragen der
Differentialgeometrie. J. Springer, Berlin, 1938.

\bibitem[B-55]{B-55} Blaschke, W.  Einführung in die Geometrie der Waben. (German) Birkhäuser Verlag, Basel-Stuttgart, 1955.

\bibitem[GS-24]{GS-24} Graf, H., Sauer, R.  \"Uber dreifache
Geradensysteme in der Ebene, welche Dreiecksnetze bilden, 
Sitzungsb. Math.-Naturw. Abt. (1924), 119-156.


\bibitem[E-74]{E-74} Erdo\v{g}an, H. I. D\"uzlemde $6$-gen doku te\c{s}kil eden \v{c}ember demety 3-\"uzleri, Ph.D. Thesis, Istanbul Teknik Ueniversitesi, Istanbul, 1974.



\bibitem[L-77]{L-77} Lazareva, V. B. Three-webs formed by families of circles on the plane.  Differential geometry (Russian), pp. 49–64, 140, Kalinin. Gos. Univ., Kalinin, 1977.

\bibitem[N-14]{N-14} Nilov, F. K. 
On new constructions in the Blaschke-Bol problem. (Russian) 
Mat. Sb. 205 (2014), no. 11, 125–144; translation in 
Sb. Math. 205 (2014), no. 11-12, 1650–1667. 

\bibitem[S-05]{S-05} Shelekhov, A. M.  Classifications of regular 3-webs formed by pencils of circles. (Russian) Sovrem. Mat. Prilozh. No. 32, Geom. Geom. Tkan. (2005), 7--28; translation in J. Math. Sci. (N.Y.) 143 (2007), no. 6, 3607--3629.

\bibitem[S-32]{S-32} Strubecker, K. 
Über eine Klasse spezieller Dreiecksnetze aus Kreisen. (German) 
Monatsh. Math. Phys. 39 (1932), no. 1, 395--398. 

\bibitem[V-29]{V-29} Volk, O.
Über spezielle Kreisnetze. (German) 
Sitzungsberichte München 59, 125--134 (1929).

\bibitem[W-38]{W-38} Wunderlich, W. \"Uber ein besonderes Dreiecksnetz aus Kreisen, Sitzungsber. Akad. Wiss. Wien 147 (1938), 385--399 


\end{thebibliography}
\end{document}